\documentstyle[11pt,amsbsy,amscd,amsfonts,amssymb,amstex,amstext]{amsart}
\theoremstyle{plain}
   \newtheorem{thm}{Theorem}[section]
   \newtheorem{lem}[thm]{Lemma}
   
   \newtheorem{prop}[thm]{Proposition}
   \newtheorem{defn-prop}[thm]{Definition-Proposition}
   \newtheorem{cor}[thm]{Corollary}
\theoremstyle{definition}
   \newtheorem{defn}[thm]{Definition}

   \newtheorem{exmp}[thm]{Example}
   
\theoremstyle{remark}
   \newtheorem{rem}[thm]{Remark}
   \newtheorem{claim}[thm]{Claim}

\begin{document}
\title[Newton number]
      {Lower bound of Newton number}
\footnote{\it{running title.} 
                 \rm{Newton number}


               }
\date{April~24~2004} 
\address[Masako Furuya]{2-155, Makinohara, Matsudo-city, Chiba, 270-2267, 
                        Japan}
\email{HZI00611@@nifty.ne.jp}
\maketitle
\begin{center} Masako {\sc Furuya} \end{center}
\begin{abstract}
  We show a lower estimate of the Milnor number of an isolated hypersurface 
  singularity, via its Newton number. We also obtain analogous estimate of 
  the Milnor number of an isolated singularity of a similar complete 
  intersection variety. 
\end{abstract}
\section*{Introduction}
   We study the Newton number of a polyhedron in order to calculate 
the Milnor number of an isolated singularity defined by an analytic mapping. 

   Section 2 treats the Newton number of a quasi-convenient polyhedron. 

   In Section 3, we consider the Milnor number $\mu(f,0)$ of an isolated 
hypersurface singular point $ 0 \in {\bold C}^n $ defined by a function 
$ f \in {\bold C}\{z_1,\dotsb,z_n\} $. It is well known that the Milnor 
number of a critical point 0 of a semi-quasihomogeneous function $f$ is equal to 
that of its initial part ([MO]-Thm.1, [A]-Thm.3.1, [LR]-Cor.2.4). On the other hand, 
Kouchnirenko proved that $ \mu(f,0) \geq \nu(f) $ holds for any function 
$ f \in {\bold C}[[z_1,\dotsb,z_n]] $, where $\nu(f)$ is the Newton number of $f$ 
([K]). We show a lower estimate of $\mu(f,0)$  
for not necessarily semi-quasihomogeneous function $f$ as follows, 
$$ \mu(f,0) \geq \nu(g) \geq (a_1-1) \cdot \dotsb \cdot (a_n-1), $$
where $(a_1,0,\dotsb,0),\dotsb,(0,\dotsb,0,a_n)$ 
are vertices of an arbitrary $n-1$ dimensional simplex lying below $\Gamma(f)$ 
with $ a_i \geq 1 \; (i=1,\dotsb,n), $ and $g$ is a standard modification 
of $f$ to a convenient function. 
When all the $a_i$ are integers, this result follows from Kouchnirenko's Theorem 
and upper semicontinuity of $\mu$ under a deformation. Recently Tomari and I 
have given a simple proof of this theorem which is quite different from 
the one of this note ([TF]). 

   In Section 4, we mention a $\mu$-constant family of an isolated 
hypersurface singularity. 

   In Section 5, we consider the case of complete intersection singularities. 
Oka obtained a formula of the Milnor number of an isolated similar complete 
intersection singularity ([O2]) which is a generalization of Kouchnirenko's 
formula ([K]). We also obtain a lower estimate of the Milnor number of 
an isolated similar complete intersection singularity. \\
{\small {\bf Acknowledgement.} 
I would like to express my sincere gratitude to Professor Mutsuo Oka 
for his valuable advice on Lemma 1.3, Corollary 1.5 and 
many useful observations throughout this paper. 
I am also grateful to Professor Kazumasa Ikeda for giving precious 
comments concerning Lemma 1.4.}
\section{Preliminary}
{\bf Notations}. 

For a subset $ I \subseteq \{1,\dotsb,n \} $, we put 
$|I| := $ the cardinality of $I$, \ $ I^c:=\{1,\dotsb,n \} \smallsetminus I $, 

$ {\bold R}^I := \{ {\bold x}=(x_1,\dotsb,x_n) \in {\bold R}^n \; ; \; 
                       x_i=0 \; $ if $ i \notin I \}, $

$ {\bold R}_I := \{ {\bold x}=(x_1,\dotsb,x_n) \in {\bold R}^n \; ; \; 
                       x_i=0 \; $ if $ i \in I \}, $

$ \pi^I \; : \; {\bold R}^n \longrightarrow {\bold R}^I, \quad 
  \pi_I \; : \; {\bold R}^n \longrightarrow {\bold R}_I $ : 
the projection maps.

For a polyhedron $X$ in $ {\bold R}^n, \; $ let \ 
$X^I:=X \cap {\bold R}^I$, 

$ V_k(X) := $ the $k$-dimensional volume of $X$ \ 
($ \dim_{\bold R}X \leq k \leq n $), 

$ \text{Vert}(X):= $ the set of all vertices of the polyhedron $X$, 

$ \text{Cone}_P(X):= $ the cone over the polyhedron $X$ with cone vertex 
$P$, \; ($ \text{Cone}_P(\emptyset):=P $). 

Let $ P_1,\dotsb,P_{r+1} $ be points in a general position in $ {\bold R}^n \; $ 
$(r \leq n)$.

$ | P_1,\dotsb,P_{r+1} | := $ 
  the r dimensional simplex generated by the r+1 vertices $ P_1,\dotsb,P_{r+1}. $

$ D^n := \{ {\bold x}=(x_1,\dotsb,x_n) \in {\bold R}^n \; ; \; 
            x_1^2 + \dotsb + x_n^2 \leq 1 \}. $ 

$ \cong \; : \; $ topologically equivalent as subspaces of Euclidean space. 
  
${\bold C} \{ {\bold z} \}:={\bold C} \{ z_1,\dotsb,z_n \}$ and 
$ {\bold z^{\lambda}}:=z_1^{\lambda_1} \cdot \dotsb \cdot z_n^{\lambda_n} $ \ 
$ ({\bold \lambda}=(\lambda_1,\dotsb,\lambda_n)). $ 
\\
\\
Recall that an arbitrary compact polyhedron can be decomposed into finite union 
of compact convex polytopes. The dimension of a convex polytope in ${\bold R}^n$ 
is defined by the dimension of the smallest $\bold R$-affine subspace of 
$ {\bold R}^n $ containing it. The dimension of a polyhedron is defined to be 
the largest possible dimension of a convex polytope contained in it. By 
convention, the dimension of the empty set is taken to be $-1$. A polyhedron 
is called pure (or pure $n$ dimensional) if it is a (not necessarily connected) 
finite union of $n$ dimensional compact convex polytopes. 
An arbitrary $n$ dimensional compact polyhedron $X$ in 
${\bold R}^n$ can be subdivided into $n$ or less than $n$ dimensional 
simplices as $ \displaystyle X = \bigcup_{t \in S}\Delta_t $ such that 
$\text{Vert}(\Delta_t) \subset \text{Vert}(X)$. 
In general, a way of subdividing the polyhedron $X$ may not be unique. 
\begin{defn}[Kouchnirenko {\rm [K]-Def.1.7}]
Let $X$ be an $n$-dimensional compact polyhedron in $ {{\bold R}_{\geq 0}}^n $. 
Then {\em the Newton number} $\nu(X)$ of $X$ is defined by 
$$ \nu(X) = \sum_{I \subseteq \{1,\dotsb,n \}} (-1)^{n-|I|} |I|! V_{|I|}(X^I).$$ 
Note that 
$ V_0(X^{\emptyset}) = 1 $ (if $ O \in X $) or 
$ V_0(X^{\emptyset}) = 0 $ (if $ O \notin X $).
\end{defn}
\begin{defn}
Assume that $ P_1,\dotsb,P_{n+1} $ are points in $({\bold R}_{\geq 0})^n$ 
in a general position and let $ X=|P_1,\dotsb,P_{n+1}| $ be the $n$ dimensional 
simplex with vertices $ P_1,\dotsb,P_{n+1} $. 
A coordinate subspace ${\bold R}^I$ of ${\bold R}^n$ is called 
{\em a full-supporting coordinate subspace} for $X$ if $ \dim_{\bold R} X^I=|I|$. 
Note that this is equivalent that ${\bold R}^I$ contains exactly 
$|I|+1$ vertices of $X$. 
\end{defn}
\begin{lem} 
Assume that $X$ is as above. Suppose that we have two full-supporting 
coordinate subspaces ${\bf R}^I$ and ${\bf R}^J$. 
Then ${\bf R}^{I\cap J}$ is also a full-supporting coordinate subspace 
for $X$. In particular, there exists a unique minimal full-supporting 
coordinate subspace ${\bf R}^{I_0}$ for $X$. 
\end{lem}
\begin{pf}
Suppose that we have two full-supporting coordinate subspaces 
${\bf R}^I$ and ${\bf R}^J$. Then it holds that $\sharp$Vert$(X^I)=|I|+1$, 
$\sharp$Vert$(X^J)=|J|+1$. Put $K:=I \cap J$ and assume that $s\le |K|+1$ 
vertices of $X$ are in ${\bf R}^I\cap {\bf R}^J={\bf R}^K$. Then we have 
$|I|+|J|+2-s$ vertices of $X$ in the coordinate subspace 
${\bf R}^{I \cup J}$. Consider the inequality \ 
$ |I|+|J|+2-s = (|I|+|J|-|K|+1)+(|K|-s+1) \ge |I|+|J|-|K|+1. $ 
By the assumption of general position, this implies that 
$s=|K|+1$ and thus ${\bf R}^{I\cap J}$ is also full-supporting for $X$.
\end{pf}
\begin{lem}
Assume that $X$ is as above such that $O \notin X$ and that ${\bf R}^I$ is 
the minimal full-supporting coordinate subspace for $X$. 
Then \ 
$ \nu(X) = |I|! V_{|I|}(X^I) \nu(\pi_I(X)). $
\end{lem}
\begin{pf} \ 
We may assume that \ $ I=\{m+1,\dotsb,n \}, \; $ 
$ X^I=|P_{m+1},\dotsb,P_{n+1}|. $ \\
Let  
$ P_i=(p_{i 1},\dotsb,p_{i n}), \; $ 
$ Q_i=(q_{i 1},\dotsb,q_{i n})=\overrightarrow{P_{n+1}P_i}, $
$$ Q=(q_{i j})_{1 \leq i \leq n, \; 1 \leq j \leq n}= 
  \begin{pmatrix}
    p_{11} & \dotsb & p_{1m} & q_{1 \; m+1}   & \dotsb & q_{1n}       \\
    \vdots &        & \vdots & \vdots         &        & \vdots       \\
    p_{m1} & \dotsb & p_{mm} & q_{m \; m+1}   & \dotsb & q_{mn}       \\
    0      & \dotsb & 0      & q_{m+1 \; m+1} & \dotsb & q_{m+1 \; n} \\
    \vdots &        & \vdots & \vdots         &        & \vdots       \\
    0      & \dotsb & 0      & q_{n \;m+1}    & \dotsb & q_{nn}
  \end{pmatrix}. $$ 
By calculation of the determinant or minor determinants of the matrix $Q$, 
we have 
$$ |J|!V_{|J|}(X^J) = 
   |I|!V_{|I|}(X^I)|J \smallsetminus I|!V_{|J \smallsetminus I|}(\pi_I(X)^J) \quad 
(I \subseteq J \subseteq \{1,\dotsb,n\}). $$
Therefore \ 
$ \displaystyle 
   \nu(X) = \sum_{I \subseteq J \subseteq \{1,\dotsb,n\}} 
            (-1)^{n-|J|} |J|! V_{|J|}(X^J) 
          = |I|! V_{|I|}(X^I) \nu(\pi_I(X)). $ 
\end{pf}
\begin{cor} \ 
Let $X$ be a pure $n$-dimensional compact polyhedron in 
$({\bold R}_{\geq 0})^n$ such that $ O \notin X $. 
Assume that there is a triangulation of X so that $X$ is a union of 
$n$-simplices $\Delta_t$ $(t\in S)$ and that 
$ {\rm Vert}(\Delta_t) \subset {\rm Vert}(X) $. 
We assume also that the simplices have the common minimal full-supporting 
coordinate subspace ${\bold R}^I$ so that 
$X \cap {\bold R}^I=\Delta_t \cap {\bold R}^I$ for any $t \in S$. 
Then we have $ \nu(X)= |I|!V_{|I|}(X^I) \nu(\pi_I(X)) $.
\end{cor}
\begin{pf} \ For a point $P$ in $\Delta_t$, we denote the minimal face 
of $\Delta_t$ containing $P$ and $X^I$ by Supp$(P, \Delta_t)$. 
\begin{claim} 
Take points $P \in \Delta_t \smallsetminus X^I$ and 
$Q \in \Delta_u \smallsetminus X^I$ and assume that $ \pi_I(P)=\pi_I(Q) $. 
Then ${\rm Supp}(P, \Delta_t)={\rm Supp}(Q, \Delta_u)$. 
In particular, 
$ \pi_I(\Delta_t) \cap \pi_I(\Delta_u)=\pi_I(\Delta_t \cap \Delta_u) $.
\end{claim}
\begin{pf}
Assume that $ \pi_I(P)=\pi_I(Q) $. 
Then it holds that $ R:=P-Q \in {\bold R}^I $. 
Let $M$ be the barycenter of $X^I$. 
Consider a relative interior point 
$P_{\epsilon}:=\epsilon P+(1-\epsilon)M $ of Supp$(P, \Delta_t)$. 
Then $P_{\epsilon}$ is also in the relative interior of 
Supp$(Q, \Delta_u)$ for a sufficiently small $\epsilon>0$ because 
$$ \epsilon P+(1-\epsilon)M
=\epsilon Q+(1-\epsilon)\left(\frac{\epsilon}{1-\epsilon}R+M \right).$$ 
This implies that ${\rm Supp}(P, \Delta_t)={\rm Supp}(Q, \Delta_u)$. 
Thus $\pi_I(P) \in \pi_I(\Delta_t \cap \Delta_u)$.
\end{pf}
Therefore $\{ \pi_I(\Delta_t) \; ; \; t \in S \}$ gives a triangulation 
of $\pi_I(X)$, and for any $ J \supseteq I $, we have 
$ \dim_{\bold R}\Delta_t^J \cap \Delta_u^J 
 = \dim_{\bold R}(\pi_I(\Delta_t)^J \cap \pi_I(\Delta_u)^J)+|I| $. 
Put $J'=J \smallsetminus I$. Thus
{\allowdisplaybreaks
\begin{align*}
\nu(X)=&\sum_{I \subseteq J \subseteq \{1,\dotsb,n\}}
                     (-1)^{n-|J|}|J|!V_{|J|}(X^J) \\
      =&\sum_{I \subseteq J \subseteq \{1,\dotsb,n\}}
                     (-1)^{n-|J|}\sum_{t \in S}\frac{1}{m_t^J}
                     |J|!V_{|J|}(\Delta_t^J) \\
      =&\sum_{I \subseteq J \subseteq \{1,\dotsb,n\}}
        (-1)^{|I^c|-|J'|}\sum_{t \in S}\frac{1}{m_t^J}
        |I|!V_{|I|}(X^I)|J'|!V_{|J'|}
        (\pi_I(\Delta_t)^{J'}) \\
      =&|I|!V_{|I|}(X^I)\sum_{I \subseteq J \subseteq \{1,\dotsb,n\}}
        (-1)^{|I^c|-|J'|}
        |J'|!V_{|J'|}(\pi_I(X)^{J'}) \\
      =&|I|!V_{|I|}(X^I)\nu(\pi_I(X)), 
\end{align*}}
where \ $ m_t^J:=\sharp\{ u \in S \; ; \; \Delta_t^J=\Delta_u^J \}$. 
\end{pf} 
\section{Newton number of quasi-convenient polyhedron}
We first note that: 
\begin{lem}
Let $X$, $Y$ be pure $n$ dimensional compact polyhedra in ${\bold R}^n$ 
with $ Y \subsetneq X $. Then the polyhedron 
$\overline{X \smallsetminus Y}$ is also pure $n$ dimensional.
\end{lem}
\begin{defn}
A polyhedron $X$ in $({\bold R}_{\geq 0})^n$ is called 
{\em quasi-convenient} if the following two conditions 
(i) and (ii) are satisfied:\\
(i) $O \in X$ and for any other vertex $P=(p_1,\dotsb,p_n)$ of $X$, 
the inequality $p_i \geq 1$ is satisfied if $p_i \ne 0$, \\
(ii) $X^I$ is topologically equivalent to $D^{|I|}$ 
for each $I \subseteq \{1,\dotsb,n \}$. \\
Note that the vertices of $X$ need not be rational. 
\end{defn}
\begin{thm}
Let $X$ be a quasi-convenient polyhedron in $({\bold R}_{\geq 0})^n$ 
and assume that $X \supseteq Y_a$ for some $a \geq 1$ where
$ Y_a:=\{{\bold x}=(x_1,\dotsb,x_n) \in ({\bold R}_{\geq 0})^n \; ; \; 
x_1+\dotsb+x_n \leq a\}. $ 
Then $ \nu(X) \geq \nu(Y_a)=(a-1)^n. $ 
\end{thm}
\begin{pf} \ 
Let $Y=Y_a$. We may assume that $\overline{X \smallsetminus Y}$ is 
not empty. So it is pure $n$ dimensional. We show the assertion 
by induction on the dimension of the polyhedron $X$. 
It is obvious for $n=1$. 
So we assume $ n \geq 2 $. Note that it holds that 
$ \nu(X)=\nu(Y)+\nu(\overline{X \smallsetminus Y}). $ \\
First, we would like to subdivide the $n$ dimensional polyhedron 
$\overline{X \smallsetminus Y}$ into pure $n$ dimensional polyhedra as 
$ \displaystyle \overline{X \smallsetminus Y}
=\bigcup_{\bar{t} \in \bar{S}} Z_{\bar{t}} $ satisfying the equality 
$$ \nu \left(\bigcup_{\bar{t} \in \bar{S}} Z_{\bar{t}} \right) 
          = \sum_{\bar{t} \in \bar{S}} \nu(Z_{\bar{t}}). $$
We construct such $n$ dimensional polyhedra $Z_{\bar{t}}$ as below. 
The polyhedron $\overline{X \smallsetminus Y}$ is subdivided into 
$n$ dimensional simplices as 
$\displaystyle \overline{X \smallsetminus Y}=\bigcup_{t \in S} \Delta_t$ 
satisfying 
$ \dim_{\bold R}\Delta_t \cap \Delta_u<n $ $(t \ne u)$ and 
$\text{Vert}(\Delta_t) \subset \text{Vert}(\overline{X \smallsetminus Y})$. 
According to Lemma 1.3, the index set $S$ can be decomposed 
as a disjoint union as follows:
$$ S = \coprod_{\emptyset \ne I \subseteq \{1,\dotsb,n\}} S(I) 
\qquad \text{where} $$ 
$$ S(I):= \{ t \in S \; ; \; {\bf R}^I \text{ is the minimal 
full-supporting coordinate subspace for } \Delta_t \}. $$
We define an equivalence relation between two elements $t$ and $u$ in 
$S$ as follows: 
$$ t \sim u \; \overset{\mathrm{def}}{\rightleftarrows} \; 
t, \; u \in S(I) \; \text{and} \; \Delta_t^I=\Delta_u^I \; 
\text{for some} \; I, $$ 
and denote the equivalence class of $t$ by $\bar{t}$. 
Then it holds  
$ \dim_{\bold R}\Delta_t^J \cap \Delta_u^J < |J| \; $ 
for any $J \subseteq \{1,\dotsb,n\}$ and $ t, \; u \in S $ such that 
$ t \not\sim u $. Thus we define \ 
$$ Z_{\bar{t}}:= \bigcup_{u \sim t}\Delta_u, \qquad  
\bar{S}:=S/\sim, $$
so that the subdivision 
$ \displaystyle \overline{X \smallsetminus Y}
=\bigcup_{\bar{t} \in \bar{S}}Z_{\bar{t}} $ 
satisfies the equality 
$ \displaystyle \nu(\overline{X \smallsetminus Y})
=\sum_{\bar{t} \in \bar{S}}\nu(Z_{\bar{t}}). $ \\
Next, we would like to show $ \nu(Z_{\bar{t}}) \geq 0 \; 
(\bar{t} \in \bar{S}). $ 
For $ t \in S(I) $, it follows from Corollary 1.5 that 
$\nu(Z_{\bar{t}})=|I|!V_{|I|}(Z_{\bar{t}}^I) \nu(\pi_I(Z_{\bar{t}}))$. 
On the other hand, the $n-|I|$ dimensional polyhedron $\pi_I(Z_{\bar{t}})$ 
is quasi-convenient. In fact, it holds for a relative interior point $P$ in 
$\Delta_t^I$ and for a sufficiently small $\epsilon>0$ that 
$ Z_{\bar{t}} \cap \overline{U_{\epsilon}(P)} = 
  ({\bold R}_{\geq 0})^n \cap \overline{U_{\epsilon}(P)} $, and thus 
$$ \pi_I(Z_{\bar{t}}) \cap \overline{U'_{\epsilon}(O)} = 
   \pi_I(Z_{\bar{t}} \cap \overline{U_{\epsilon}(P)}) = 
   \pi_I(({\bold R}_{\geq 0})^n \cap \overline{U_{\epsilon}(P)}) = 
   ({{\bold R}_{\geq 0}}^n)^{I^c} \cap \overline{U'_{\epsilon}(O)} 
   \cong D^{n-|I|}, $$ 
where $U_{\epsilon}(P)$ (resp. $U'_{\epsilon}(O)$) is the open 
$\epsilon$-neighbourhood of $P$ in ${\bold R}^n$ 
(resp. of $O$ in ${\bold R}_I)$. \\
Hence 
$ \pi_I(Z_{\bar{t}}) = \bigcup_{u \sim t} \pi_I(\Delta_u) 
  \cong D^{n-|I|}. $ 
Similarly we have 
$ \pi_I(Z_{\bar{t}})^J \cong D^{|J|} $ for $ J \subseteq I^c $. 
Therefore, by inductive hypothesis, it holds 
$ \nu(\pi_I(Z_{\bar{t}})) \geq 0, $ and so $ \nu(Z_{\bar{t}}) \geq 0. $  
\end{pf}
\begin{cor} 
We have the following assertion for a quasi-convenient polyhedron $X$ 
in $({\bold R}_{\geq 0})^n$. 
(i) If $\nu(X)=0$ then 
$E_j=(0,\dotsb,0,\underset{j}{\underset{\frown}{1}},0,\dotsb,0) \in$ 
{\rm Vert}$(X)$ for some $j \in \{1,\dotsb,n\}$. \\ 
(ii) Then the equality $ \nu(X)=0 $ holds if $E_j \in$ {\rm Vert}$(X)$ for 
some $j \in \{1,\dotsb,n\}$ and the other vertices are in ${\bold R}_{\{j\}}$.\\
(iii) Assume further that $({\bold R}_{\geq 0})^n \smallsetminus X$ is convex. 
Then the equality $\nu(X)=0$ holds if and only if 
$E_j \in $ {\rm Vert}$(X)$ for some $j \in \{1,\dotsb,n\}$. 
\end{cor}
\begin{pf}
First we prove the assertion (i). Suppose that $E_j \notin $ Vert$(X)$ 
for any $ j \in \{1,\dotsb,n\} $. Then the polyhedron $X$ includes 
a subset $Y_{1+\epsilon}$ for a sufficiently small $\epsilon>0$. 
This proves the assertion (i). 
We show the assertion (ii) by induction on the dimension of $X$. 
Assume $E_j \in$ Vert$(X)$ and 
Vert$(X) \smallsetminus \{E_j\} \subset {\bold R}_{\{j\}}$. 
Let $Y=Y_1$. Note that $Y_1 \subseteq X$. For any subdivision 
$ \displaystyle X=Y\cup\bigcup_{t \in S}\Delta_t $ 
as in Proof of Theorem 2.3, each simplex $\Delta_t$ has vertices in 
${\bold R}_{\{j\}}$ except $E_j$. Thus the minimal full-supporting coordinate 
subspace ${\bold R}^I$ satisfies $ j \notin I $. Take an index $t \in S(I)$. 
By Corollary 1.5, we get 
$\nu(Z_{\bar{t}})=|I|!V_{|I|}(Z_{\bar{t}}^I) \nu(\pi_I(Z_{\bar{t}}))$. 
Note that $E_j$ is again a vertex of $\pi_I(\Delta_t)$ and the other vertices 
of $\pi_I(\Delta_t)$ are in ${\bold R}^{I^c \smallsetminus \{j\}}$. 
Thus by the inductive hypothesis, we have $\nu(Z_{\bar{t}})=0$. 
We consider the assertion (iii). Assume that $E_j \in $ Vert$(X)$ and 
$({\bold R}_{\geq 0})^n \smallsetminus X$ is convex. Then we assert that 
Vert$(X) \smallsetminus \{E_j\} \subset {\bold R}_{\{j\}}$. 
In fact, suppose that Vert$(X) \smallsetminus \{E_j\} 
\not\subset {\bold R}_{\{j\}}$ 
and that the j-th coordinate of a point 
$P \in $ Vert$(X) \smallsetminus \{E_j\}$ 
is equal to the maximum of the j-th coordinates of all vertices of $X$. 
Then $({\bold R}_{\geq 0})^n \smallsetminus X$ 
is non-convex around the point $P$.
\end{pf}
The following Corollary is a generalization of Theorem 2.3. 
\begin{cor} 
Let $ X, \; Y $ be quasi-convenient polyhedra in $({\bold R}_{\geq 0})^n$ 
with $ Y \subseteq X.$ \\ Then $ \nu(X) \geq \nu(Y) \geq 0 $ and 
$ \nu(\overline{X \smallsetminus Y}) \geq 0 $. 
In particular, if $Y=|O,A_1,\dotsb,A_n|$, \\
$A_i=(0,\dotsb,0, \underset{i}{\underset{\frown}{a_i}}, 0,\dotsb,0)$ 
$ (1 \leq i \leq n) $, then 
$ \nu(X) \geq (a_1-1)\cdot\dotsb\cdot(a_n-1) \geq 0. $ 
\end{cor}
\begin{pf}
Put $Z:=\overline{X \smallsetminus Y}$. Then it holds that 
$\dim_{\bold R}Y^I \cap Z^I<|I|$ for each $I \subseteq \{1,\dotsb,n\}$. 
Therefore $\nu(X)=\nu(Y)+\nu(Z)$. 
The rest is similar to Proof of Theorem 2.3.
\end{pf}
\section{Lower bound of Milnor number of hypersurface singularity}
As a result of Kouchnirenko's formula ([K]-Thm.I) and Corollary 2.5 above, 
we obtain a lower estimate of Milnor number of an isolated hypersurface 
singularity (Corollary 3.1). Using the same notation as in [K], 
let $\Gamma(f)$ be the Newton boundary of $f$ with respect to 
the fixed coordinates $(z_1,\dotsb,z_n)$, 
and let $\Gamma_-(f)$ be the cone of $\Gamma(f)$ and the origin. 
An analytic function $f$ is called 
{\em convenient} if the intersection of $\Gamma(f)$ and 
each coordinate axis is non-empty. Note that if a function $f$ 
is convenient then $\Gamma_-(f)$ is quasi-convenient. 
\begin{cor}
Let $ f \in {\bold C} \{ z_1, \dotsb, z_n \} $ define an isolated 
singularity at the origin $ 0 \in {\bold C}^n $. 
Let $H$ be an arbitrary hyperplane lying below $\Gamma(f)$ 
and intersecting coordinate axis $ x_1, \dotsb, x_n $ 
at points $ a_1, \dotsb, a_n, $ 
where $ a_i \geq 1 \; (i=1,\dotsb,n) $. Then 
$ \mu(f,0) \geq \left( a_1-1 \right) \cdot \dotsb
           \cdot \left( a_n-1 \right). $
\end{cor}
\begin{pf}  
Take a sufficiently large positive integer $m$ and let 
$g:=f+z_1^m+\dotsb+z_n^m$ be a standard modification of $f$ 
to a convenient function. Then it follows from Kouchnirenko's formula 
([K]-Thm.I) and Corollary 2.5 that 
$ \mu(f,0) = \mu(g,0) \geq \nu(\Gamma_-(g)) 
  \geq (a_1-1) \cdot \dotsb \cdot (a_n-1). $ 
\end{pf}
\begin{rem}
The main results Theorem 2.3, Corollary 2.5 and 3.1 for the case 
that all the vertices of the polyhedra $ X, \; Y, \; H $ are lattice points 
follow from Kouchnirenko's Theorem ([K]-Theorem.I) and upper semicontinuity 
of $\mu$ under a deformation ([M]).
\end{rem}
\begin{rem}
Tomari has recently given a simple proof of Corollary 3.1 which uses 
a theory of multiplicity of filtered ring. He has also proved that 
the equality holds if and only if $f$ is semi-quasihomogeneous ([TF]). 
\end{rem}
\section{$\mu$-constant family of three dimensional hypersurface singularity}
Let $f({\bold z})=\sum_{\lambda \in \Lambda}\gamma_{\lambda}{\bold z}^{\lambda} 
\in {\bold C}\{{\bold z}\}$ be an analytic function. 
As an application of Corollary 2.4 and 2.5 above, we investigate 
the family of {\em a negligible truncation} defined in [O1]. 

Let $A$ be a vertex of $\Gamma(f)$ and put 
$f_t({\bold z})=f({\bold z})-(1-t)\gamma_A{\bold z}^A$. We assume that both 
$f_1$ and $f_0$ are convenient and that $\Gamma_-(f_1)$ is a proper subset 
of $\Gamma_-(f_0)$. Then we see that 
$\overline{\Gamma_-(f_0) \smallsetminus \Gamma_-(f_1)}
=\text{Cone}_A(\overline{\Gamma(f_0) \smallsetminus \Gamma(f_1)})$. 
Therefore if $A \notin {\bold R}^I$ then 
$\overline{\Gamma_-(f_0)^I \smallsetminus \Gamma_-(f_1)^I}$ is 
strictly less than $|I|$ dimensional. If $A \in {\bold R}^I$ then 
$\overline{\Gamma_-(f_0)^I \smallsetminus \Gamma_-(f_1)^I}$ 
is pure $|I|$ dimensional because 
$\Gamma_-(f_0)^I \ne \Gamma_-(f_1)^I$ and both 
$\Gamma_-(f_0)^I$ and $\Gamma_-(f_1)^I$ are pure $|I|$ dimensional. 

Now we consider the case that $n=4$ and 
${\bold C}\{{\bold z}\}={\bold C}\{x,y,z,w\}$. 
By account in the preceding paragraph, we have the following Proposition 
as a particular case of Corollary 2.5. 
\begin{prop}
Let $A=(a_1,a_2,a_3,a_4)$ and assume that 
$\Gamma_-(f_0)=\Gamma_-(f_1) \cup \Delta$ holds for some 
$4$-simplex $\Delta=|A, B, C, D, E|$. 
Then the equality $\nu(\Gamma_-(f_0))=\nu(\Gamma_-(f_1))$ holds 
if and only if one of the following (i)-(iii) is satisfied 
(up to a permutation of the coordinates). \\
(i) Case $a_1=0$, $a_i>0$ $(i=2,3,4)$. 
Then $E=(1,*,*,*)$.\\
(ii) Case $a_1=a_2=0$, $a_3>0$, $a_4>0$. 
Then either $D=(0,1,*,*)$ or $E=(1,0,*,*)$. \\
(iii) Case $a_1=a_2=a_3=0$, $a_4>0$. 
Then either $C=(0,0,1,*)$ or $D=(0,1,0,*)$ or $E=(1,0,0,*)$.
\end{prop}
\begin{pf} \ 
If $a_i>0$ $(i=1,2,3,4)$ then ${\bold R}^4$ is the minimal 
full-supporting coordinate subspace for $\Delta$, hence $\nu(\Delta)>0$. 
Thus the vertex $A$ must be as in one of (i)-(iii) above.\\
Case (i). Then ${\bold R}^{\{2,3,4\}}$ is the minimal full-supporting 
coordinate subspace for $\Delta$. The equality $\nu(\Delta)=0$ holds 
if and only if $\pi_{\{2,3,4\}}(\Delta)=|O,E_1|$ where 
$E_1:=(1,0,0,0)$. \\
Case (ii). Then ${\bold R}^{\{3,4\}}$ is the minimal full-supporting 
coordinate subspace for $\Delta$. The equality $\nu(\Delta)=0$ holds 
if and only if the $2$-simplex $\pi_{\{3,4\}}(\Delta)$ 
has either $E_1$ or $E_2:=(0,1,0,0)$ as a vertex. \\
Case (iii). Then ${\bold R}^{\{4\}}$ is the minimal full-supporting 
coordinate subspace for $\Delta$. The equality $\nu(\Delta)=0$ holds 
if and only if the $3$-simplex $\pi_{\{4\}}(\Delta)$ 
has at least one of $E_1$ or $E_2$ or $E_3:=(0,0,1,0)$ as a vertex. 
\end{pf}
Each family $f_t$ in the following examples is negligible truncation, 
so it is $\mu$-constant.
\begin{exmp} 
$f_t=x^3+y^3+z^5+xw^5+ty^2zw+w^m$ \ \ $(m \geq 8)$. 
\end{exmp}
\begin{pf}
It holds that $\Gamma_-(f_0)=\Gamma_-(f_1) \cup \Delta$ for the 
$4$-simplex $\Delta$ with the vertices $(0,3,0,0)$, $(0,0,5,0)$, 
$(1,0,0,5)$, $(0,2,1,1)$, $(0,0,0,m)$, and ${\bold R}^{\{2,3,4\}}$ 
is the minimal full-supporting \\coordinate subspace for $\Delta$. 
We see that $\pi_{\{2,3,4\}}(\Delta)=|O,E_1|$ and so $\nu(\Delta)=0$. 
\end{pf}
\begin{exmp}
$f_t=x^3+y^3+z^5+xw^5+yw^5+tzw^6+w^m$ \ \ $(m \geq 8)$. 
\end{exmp}
\begin{pf}
It holds that $\Gamma_-(f_0)=\Gamma_-(f_1) \cup \Delta$ for the 
$4$-simplex $\Delta$ with the vertices $(0,0,5,0)$, $(1,0,0,5)$, 
$(0,1,0,5)$, $(0,0,1,6)$, $(0,0,0,m)$, and ${\bold R}^{\{3,4\}}$ is the 
minimal full-supporting coordinate subspace for $\Delta$. We see that 
$\pi_{\{3,4\}}(\Delta)=|O,E_1,E_2|$ and so $\nu(\Delta)=0$. 
\end{pf}
\begin{exmp} 
$f_t=x^2+y^5+z^6+yw^6+z^2w^5+tw^8+w^m$ \ \ $(m \geq 9)$. 
\end{exmp}
\begin{pf}
$\Gamma_-(f_0)=\Gamma_-(f_1) \cup \Delta$ for 
$\Delta=|(2,0,0,0), (0,1,0,6), (0,0,2,5), (0,0,0,8), (0,0,0,m)|$. 
Therefore 
$\pi_{\{4\}}(\Delta)=|(0,0,0,0), \; (2,0,0,0), \; (0,1,0,0), \; (0,0,2,0)|$ 
and $\nu(\Delta)=0$.
\end{pf}
\section{The r-th Newton number and similar complete intersection singularity}
The r-th Newton number is a natural generalization of the Newton number 
([O2]-Thm.7.2). Let $d_1, \dotsb, d_r$ be positive integers. 
For an $n$-dimensional compact polyhedron $X$ in 
${{\bold R}_{\geq 0}}^n$, 
{\em the r-th Newton number of $(d_1, \dotsb, d_r)$-type} 
${\nu}^r_{d_1 \dotsb d_r}(X)$ of $X$ \ 
$(1 \leq r \leq n)$ is defined by 
$$ {\nu}^r_{d_1 \dotsb d_r}(X) =  
\sum_{I \subseteq \{1,\dotsb,n \}, \; r \leq |I|}(-1)^{n-|I|} 
F^{|I|}_r(d_1,\dotsb,d_r)|I|! V_{|I|}(X^I) + \epsilon(-1)^{n-r+1}, $$ 
where $ \epsilon=1 $ (if $ O \in X $) or $ \epsilon=0 $ (if $ O \notin X $), 
$ \displaystyle 
F^l_k(d_1,\dotsb,d_k):=\sum_{i_1+\dotsb+i_k=l-k}{d_1}^{i_1+1}\dotsb{d_k}^{i_k+1}.$ 
Then the following is shown in the same way as the case of the Newton number. 
\begin{thm}
Let $ X \subset ({\bold R}_{\geq 0})^n $ be quasi-convenient. 
Then $ {\nu}^r_{d_1 \dotsb d_r}(X) \geq 0 $ for $ 1 \leq r \leq n $. 
\end{thm}
\begin{cor}
Let $ X, \; Y $ be quasi-convenient polyhedra in $({\bold R}_{\geq 0})^n$ 
with $ Y \subseteq X.$ \\
Then $ {\nu}^r_{d_1 \dotsb d_r}(X) \geq {\nu}^r_{d_1 \dotsb d_r}(Y) \geq 0 $ 
and  $ {\nu}^r_{d_1 \dotsb d_r}(\overline{X \smallsetminus Y}) \geq 0 $ 
for $ 1 \leq r \leq n $. \\
In particular, if $Y=|O,A_1,\dotsb,A_n|$, 
$ A_i=(0,\dotsb,0, \underset{i}{\underset{\frown}{a_i}}, 0,\dotsb,0)$ 
$ (i=1,\dotsb,n) $, then 
$$ {\nu}^r_{d_1 \dotsb d_r}(X) \geq 
   \sum_{s=r}^n (-1)^{n-s} F^s_r(d_1,\dotsb,d_r) \sigma_s(a_1,\dotsb,a_n)
   +(-1)^{n-r+1} \geq 0, $$
where $ \sigma_s(a_1,\dotsb,a_n) $ is the s-th elementary symmetric 
function of $ a_1,\dotsb,a_n $ defined by 
$$ \prod_{i=1}^n(t-a_i)=t^n-\sigma_1 t^{n-1}+\dotsb+(-1)^n \sigma_n. $$
\end{cor}
The key Lemma for Proof of Theorem 5.1 is the following 5.3 which is 
an application of Corollary 1.5. Using this Lemma, we can show 
Theorem 5.1 by the same arguments of Proof of Theorem 2.3. 
Let $ \displaystyle 
G^l_k(d_1,\dotsb,d_k):=\sum_{i_1+\dotsb+i_k=l-k}{d_1}^{i_1}\dotsb{d_k}^{i_k}. $ 
\begin{lem}
Let $X$, $I$ be as in Corollary 1.5. Set $m=n-|I|$, $X'=\pi_I(X)$. 
Then the following hold for $ 1 < r < n $. \ 
If $ r \leq |I| $ then 
{\allowdisplaybreaks
\begin{align*}
{\nu}^r_{d_1 \dotsb d_r}(X)= 
              & \sum_{I \subseteq J \subseteq \{1,\dotsb,n \}}(-1)^{n-|J|} 
                F^{|J|}_r(d_1,\dotsb,d_r)|J|! V_{|J|}(X^J) \\
            = & |I|!V_{|I|}(X^I) \left\{ 
               {d_r}^{n-m}{\nu}^r_{d_1 \dotsb d_r}(X')
              +d_r G^{n-m+1}_2(d_{r-1}, \; d_r){\nu}^{r-1}_{d_1 \dotsb d_{r-1}}(X')
              \right. \\
              & \left. \hspace*{2.2cm} +\dotsb 
              +d_r \dotsb d_2 G^{n-m+1}_r(d_1, \dotsb ,d_r){\nu}^1_{d_1}(X') 
              \right\} \quad (r \leq m), \; \; \text{or} \\
              & |I|!V_{|I|}(X^I) \left\{ 
               d_r \dotsb d_{m+1}G^{n-m+1}_{r-m+1}(d_m, \dotsb, d_r)
                                           {\nu}^m_{d_1 \dotsb d_m}(X')
              +\dotsb \right. \\
              & \left. \quad 
              +d_r \dotsb d_2 G^{n-m+1}_r(d_1, \dotsb, d_r){\nu}^1_{d_1}(X')
              +F^{n-m}_{r-m}(d_{m+1}, \dotsb, d_r) \right\} \quad (r>m). 
\end{align*}}
If $ r>|I| $ then 
{\allowdisplaybreaks
\begin{align*}
{\nu}^r_{d_1 \dotsb d_r}(X)= 
              & \sum_{I \subseteq J \subseteq \{1,\dotsb,n \}, \; r \leq |J|}(-1)^{n-|J|} 
                F^{|J|}_r(d_1,\dotsb,d_r)|J|! V_{|J|}(X^J) \\
            = & |I|!V_{|I|}(X^I) \left\{ 
              G^{n-m+1}_1(d_r){\nu}^r_{d_1 \dotsb d_r}(X')
              +d_r G^{n-m+1}_2(d_r,d_{r-1}){\nu}^{r-1}_{d_1 \dotsb d_{r-1}}(X')
              \right. \\
              & \hspace*{2.2cm}  
              \left. +\dotsb 
              +d_r \dotsb d_{r-n+m+1}G^{n-m+1}_{n-m+1}(d_r, \dotsb, d_{r-n+m})
                  {\nu}^{r-n+m}_{d_1 \dotsb d_{r-n+m}}(X') \right\} \\
              & \hspace*{6.3cm} \quad (r \leq m), \quad \text{or} \\
              & |I|!V_{|I|}(X^I) \left\{ 
              d_r \dotsb d_{m+1}G^{n-m+1}_{r-m+1}(d_r, \dotsb, d_m)
                            {\nu}^m_{d_1 \dotsb d_m}(X')+\dotsb \right. \\
              & \hspace*{2.2cm} 
              +d_r \dotsb d_{r-n+m+1}G^{n-m+1}_{n-m+1}(d_r, \dotsb, d_{r-n+m})
                       {\nu}^{r-n+m}_{d_1 \dotsb d_{r-n+m}}(X') \\
              & \hspace*{2.2cm} \left. 
              +F^{n-m}_{r-m}(d_{m+1}, \dotsb, d_r) \right\} \quad (r>m). 
\end{align*}}
\end{lem}
Next we would like to give a similar result as Corollary 3.1 for 
more general situation. 
Let $ {\bold f}=(f_1,\dotsb,f_r) \; : \; ({\bold C}^n,0) \; \to \; 
({\bold C}^r,0) \; (2 \leq r < n) $ be a germ of an analytic mapping 
such that $ {\bold f}(0)=0 $. We assume that ${\bold f}^{-1}(0)$ is 
a germ of {\em a similar complete intersection variety} (SCIV) 
with an isolated singularity at the origin $ 0 \in {\bold C}^n $. 
(See [O2]-7.) 
As a consequence of Corollary 5.2 and Oka's formula, 
we obtain a lower estimate of the Milnor number $\mu({\bold f},0)$. 
\begin{thm}[Oka {\rm [O2]-Thm.7.2}]
Let $ {\bold f}=(f_1,\dotsb,f_r) \; : \; ({\bold C}^n,0) \; \to \; 
({\bold C}^r,0) $ define a germ of a SCIV such that 
$ \Gamma(f_j)=d_j\Gamma(f) $ 
with an isolated singularity at the origin. Then it holds that 
$ \mu({\bold f},0) \geq {\nu}^r_{d_1 \dotsb d_r}(\Gamma_-(f)). $ 
If ${\bold f}^{-1}(0)$ is non-degenerate then the equality holds. 
\end{thm}
\begin{cor}
Let \ $\bold f$ be as in Theorem 5.4. 
Suppose that $ \Gamma_-(f) \supseteq |O,A_1,\dotsb,A_n| $, \\
$ A_i=(0,\dotsb,0, \underset{i}{\underset{\frown}{a_i}}, 0,\dotsb,0)$, 
$ a_i \geq 1 $ $(i=1,\dotsb,n)$. Then 
$$ \mu({\bold f},0) \geq 
   \sum_{s=r}^n (-1)^{n-s} F^s_r(d_1,\dotsb,d_r) \sigma_s(a_1,\dotsb,a_n)
   +(-1)^{n-r+1}. $$
Particularly, if $d_1=\dotsb=d_r=1$ then 
$$ \mu({\bold f},0) \geq 
   \sum_{s=r}^n (-1)^{n-s} \binom{s-1}{r-1} \sigma_s(a_1,\dotsb,a_n)
   +(-1)^{n-r+1}. $$
\end{cor}
%

\end{document}